\documentclass[12pt]{article}
\usepackage[utf8]{inputenc}
\usepackage{amsmath}
\usepackage{bbm}
\usepackage{mathtools}
\usepackage{amsfonts}
\usepackage{dsfont}
\usepackage{amssymb}
\usepackage[margin=1in,headheight=1in, asymmetric]{geometry}
\usepackage{amsthm}
\usepackage{color}

\newtheorem{theorem}{Theorem}
\newtheorem{proposition}[theorem]{Proposition}
\newtheorem{lemma}[theorem]{Lemma}
\newtheorem{definition}[theorem]{Definition}
\newtheorem{example}[theorem]{Example}
\newtheorem{corollary}[theorem]{Corollary}

\newtheorem{remark}[theorem]{Remark}

\newcommand{\ee}{\mathbb{E}}

\usepackage{comment}
\usepackage{setspace}

\newcommand{\tr}{\mathrm{T}}

\title{Self-similar Gaussian Markov processes}
\author{Benedict Bauer\thanks{Financial support from the Austrian Science Fund (FWF) under grants P~30750
and Y~1235 is gratefully acknowledged.}\\
University of Vienna \\
\tt{benedict.bauer@univie.ac.at}\\
\\
 Stefan Gerhold \\
TU Wien \\
\tt{sgerhold@fam.tuwien.ac.at}
}

\date{\today}

\numberwithin{equation}{section}
\numberwithin{theorem}{section}

\begin{document}
%\onehalfspacing
%\noindent

\maketitle

\begin{abstract}
We characterize all multi-dimensional real self-similar Gaussian Markov processes.
Three types of covariance matrix functions occur: white-noise type functions, covariances
that can be expressed by continuous matrix semigroups, and covariances based on non-continuous solutions of Cauchy's functional equation. Characterizing the latter requires us
to develop some results on the representation theory of non-continuous matrix semigroups,
which are presented in a companion paper.
In dimension one, besides white noise, the self-similar Gaussian Markov processes
reduce to a two-parameter family of time-changed Brownian motions.  
This observation simplifies several proofs of non-Markovianity of concrete processes found in the literature.
\end{abstract}

\smallskip


\section{Introduction}

Fractional Brownian motion with Hurst parameter~$H$ is among the best-known examples
of $H$-self-similar Gaussian stochastic processes. It is not Markovian for $\tfrac12\neq H\in(0,1)$, which can be proven
using the functional equation that the covariance function of any Gaussian Markov 
process must satisfy (see~\cite{Hu03} and~\cite[Theorem~2.3]{No12}).
Many authors have considered self-similar Gaussian processes that are variants and extensions of fractional Brownian
motion~\cite{BoGoTa07,RuTu06,Sg13,Sg14,Zi17}, and for some, the natural question of Markovianity
has been studied.
In dimension one, it is known that all Gaussian Markov processes are time-changed
Brownian motions~\cite[p.~4]{AdCaSa90}. Our contribution in dimension one is that we explicitly
determine this time change for self-similar processes; it depends on the self-similarity parameter and
an additional parameter (Theorem~\ref{thm:1dim}). It is usually a very simple matter
to compare the resulting explicit covariance function with a given covariance function,
yielding easy non-Markovianity proofs for many self-similar Gaussian processes that
have been considered in the literature. 
Moreover, based on a preprint of the present paper, the two-parameter family
of Gaussian processes
we just mentioned has been analyzed in detail
in~\cite{El24}, where it is proposed as a model for
anomalous diffusion.

Our main contribution is that we determine all \emph{multi-dimensional} self-similar Gaussian Markov
processes.
By self-similarity, the bivariate covariance function can
be reduced to a function of a single argument, and the functional equation mentioned
above shows that this function must satisfy the semigroup property $g(x+y)=g(x)g(y)$.
We identify a contractivity estimate involving
the self-similarity parameter, which must be satisfied for self-similar Gaussian
Markov processes with continuous semigroup.
Under mild assumptions, \emph{all} matrix solutions of the equation $g(x+y)=g(x)g(y)$, including non-measurable
ones, were determined in~\cite{BaGe24}.
While the self-similar Gaussian Markov processes resulting from these non-measurable semigroups
 are not likely to occur in
applications, they have to be studied to obtain a full classification.

The paper is organized as follows.  In Section~\ref{se:prelim}, we define the class of 
stochastic processes we wish to consider, and state our main characterization results.
They are proven in Section~\ref{se:cov}, which includes a new criterion for
checking positive definiteness in the multi-dimensional case.
Section~\ref{se:volt} presents Volterra representations for self-similar Gaussian
Markov processes, assuming that the associated semigroup is continuous,
and that the contractivity estimate mentioned above is strict.
Our multidimensional characterization statements are specialized to dimension one in Section~\ref{se:1dim},
which includes a sketch of a simpler proof for this special case. 
Section~\ref{se:non-markov} returns to our initial motivation, by giving straightforward
proofs that some self-similar Gaussian processes from the literature are not Markovian.

\section{Preliminaries and main results}\label{se:prelim}

Throughout this section and the following one, $V$ denotes a real $d$-dimensional vector space
equipped with an inner product $\langle \cdot, \cdot\rangle$. We write $L(V)$  for the set of linear maps from~$V$ to itself.
The inner product $\langle \cdot, \cdot\rangle_X$ associated with a $V$-valued Gaussian process~$X$
will be defined in Proposition~\ref{prop:indep}, and we will write  $\|\cdot\|_X$ for the corresponding
 operator norm on $L(V)$. The adjoint of a linear map will be denoted by $^*$, where the inner product
 will be clear from the context, or explicitly stated otherwise.
%and $\|\cdot\|_{\mathrm{op}}$ for the operator norm.
%\inlinecomment{$\mathrm{End}(V)$ Standard-Notation?}
\begin{definition}
A stochastic process $(X_t)_{t \geq 0}$ with values in $V$
 is called Gaussian if, for any finite set $S \subset \mathbb{R}_{\ge 0}$, the $V^{|S|}$-valued random variable $(X_s)_{s \in S}$ is Gaussian.
\end{definition}
It is well known that a Gaussian process $X$ with values in $V$ is uniquely characterized in law by the two mappings $t\mapsto \ee[X_t]$
and $(s,t)\mapsto R(s,t) \in L(V)$ given by
\[
\langle v,R(s,t)u\rangle=\ee[\langle X_s,v\rangle \langle X_t,u \rangle].
\]
The mapping $R$ is called the \textit{covariance kernel} of $X$. We also write $R_X$, if~$X$ might
not be clear from the context. 
%Conversely, for any function~$R(\cdot,\cdot)$
% satisfying $R(s,t)=R(t,s)^\tr$ and the positive definiteness property %for any $n \in \mathbb{N}$
%\begin{align}
%\label{condition}
%\sum_{i,j=1}^n \langle v_i,R(t_i,t_j)v_j \rangle \ge 0
%\end{align}
%for all $v_1,\dots,v_n \in V$ and $t_1,\dots,t_n \in \mathbb{R}_{\ge 0}$, there exists a unique centered Gaussian process with covariance kernel $R$. We can characterize Markovianity for Gaussian processes via their covariance kernels.
%
\begin{definition}
For $H>0$ we call $X$ $H$-self-similar if $(X_{at})_{t \ge 0} \sim (a^{H}X_t)_{t \ge 0}$ 
(equality in distribution) for any $a>0$, and refer to~$H$ as the self-similarity parameter.
\end{definition}
It is clear, and well-known, that $H$-self-similarity of a Gaussian process
can be characterized via the covariance kernel.
\begin{theorem}
\label{similar}
Let $X$ be a centered Gaussian process. Then $X$ is $H$-self-similar if and only if $R$ satisfies
\[
R(as,at)=a^{2H}R(s,t)
\]
for any $a,s,t>0$.
\end{theorem}
\begin{corollary}
Let $X$ be an $H$-self-similar centered Gaussian process on $V$. Then $\operatorname{supp}(X_t)=\operatorname{ker}(R(1,1))^\perp$ for all $t>0$.
\end{corollary}
\begin{proof}
By self-similarity, we have $\operatorname{ker}(R(1,1))=\operatorname{ker}(R(t,t))$ for any $t>0$.
Denote by $p: V \to V$ the orthogonal projection onto $\operatorname{ker}(R(1,1))$ and set $Y_t := pX_t$.  Clearly, $Y$ is again Gaussian since the orthogonal projection is linear.  Denote by $R_Y$ the covariance kernel of~$Y$. We have
\[
\langle v_1, R_Y(s,t)v_2 \rangle = \ee[\langle v_1,pX_s\rangle \langle v_2,pX_t\rangle],
\quad v_1,v_2 \in V.
\]
Since orthogonal projections are self-adjoint, we obtain
\[
\ee[\langle v_1,pX_s\rangle \langle v_2,pX_t\rangle]=\ee[\langle pv_1,X_s\rangle \langle pv_2,X_t\rangle]=\langle v_1, pR(s,t)pv_2 \rangle,
\]
where $R$ is the covariance kernel of $X$.  Since $pv_2\in \operatorname{ker}(R(1,1))$, we obtain $R_Y(t,t)=0$ for any $t\ge 0$. This easily implies that~$Y$ is the zero process.
%By Lemma~\ref{UpperBound} we have $R_Y(s,t)=0$ for all $s,t \ge 0$ and hence by uniqueness of the covariance kernel we have $Y_t = 0$ a.s.
%\inlinecomment{Das Lemma sollte nicht noetig sein. $R_Y(t,t)=0$, also fuer jedes $t$ Varianz 0, dann ist auch fuer $s,t$ das Produkt $Y_s Y_t$ f.s. null.}
\end{proof}
If $R(t,t)=t^{2H}R(1,1)$, $t>0$, is not invertible, then the corollary shows that we are dealing with
a self-similar Gaussian
 process of lower dimension, and~$V$ can be replaced by $\operatorname{ker}(R(1,1))^\perp$.
Thus, it suffices to characterize processes with full support. In this case,
we can define the inner product
\begin{equation}\label{eq:def ip}
\langle v,u \rangle_X = \langle v, R(1,1)^{-1}u \rangle,
\end{equation}
where $\langle \cdot, \cdot \rangle$ is the original inner product on $V$. This is well-defined, since $R(1,1)^{-1} \in L(V)$ is symmetric and positive definite.
\begin{proposition}\label{prop:indep}
Let~$X$ be a $V$-valued centered Gaussian process with $R(1,1)$ invertible. Then
the inner product $\langle \cdot,\cdot \rangle_X=\langle \cdot, R(1,1)^{-1}\cdot \rangle$ is independent of the choice of $\langle \cdot, \cdot \rangle$.
\end{proposition}
\begin{proof}
Let $\langle \cdot, \cdot \rangle$ be an inner product on $V$. For any other inner product on $V$ there exists a symmetric positive map $A$ such that it is given by $\langle \cdot, A\cdot \rangle=: \langle \cdot, \cdot \rangle_A$.
%Given an SSGM process $X$, denote by $R(s,t)$ the covariance kernel with respect to the inner product $\langle \cdot, \cdot \rangle$. 
Calculating the covariance kernel with respect to $\langle \cdot, \cdot \rangle_A$, we obtain
\begin{align*}
    \ee[\langle Av,X_s \rangle \langle Au,X_t \rangle]&= \langle Av, R(s,t)Au\rangle\\
     &= \langle v, AR(s,t)Au\rangle= \langle v, R(s,t)Au\rangle_A,
     \quad u,v \in V.
\end{align*}
Hence the covariance kernel with respect to $\langle \cdot, \cdot \rangle_A$ is $R(s,t)A$. We calculate
\[
\langle v, (R(1,1)A)^{-1}u \rangle_A= \langle v, A^{-1}R(1,1)^{-1}u \rangle_A=\langle v, R(1,1)^{-1}u \rangle.
\]
Hence the definition of $\langle \cdot, \cdot \rangle_X$ does not depend on our choice of $\langle \cdot, \cdot \rangle$.
\end{proof}
In the notation of the preceding proof, the inner product~\eqref{eq:def ip} induced by~$X$
is an abbreviation for $\langle \cdot, \cdot \rangle_{R_X(1,1)^{-1}}$.
\begin{lemma}\label{le:identity}
Let~$X$ be as in Proposition~\ref{prop:indep}.
Under $\langle \cdot, \cdot \rangle_X$, the covariance kernel of~$X$ is $R(s,t)R(1,1)^{-1}$. 
\end{lemma}
\begin{proof}
Put $A=R(1,1)^{-1}$ in the proof of the preceding proposition.
\end{proof}
\begin{definition}\label{def:ssgm}
  For $H>0$, we call a $V$-valued process~$X$ an $H$-SSGM process, or just SSGM process, if
  \begin{itemize}
    \item[(i)] it is a centered $H$-self-similar Gaussian Markov process, and
    \item[(ii)] it has full support, i.e., $R(1,1)$ is invertible.
  \end{itemize}
  For such a process, it is understood that the finite-dimensional
  vector space~$V$ is equipped with the inner product $\langle \cdot,\cdot \rangle_X$
  and the corresponding operator norm.
  We write $\tilde R$, or sometimes $\tilde{R}_X$, for the covariance kernel w.r.t.\ $\langle \cdot,\cdot \rangle_X$,
  so that $\tilde{R}(1,1)=\operatorname{id}$ by Lemma~\ref{le:identity}.
\end{definition}
%
%From now on we
%often just write $\langle \cdot, \cdot \rangle$ instead of
%$\langle \cdot, \cdot \rangle_X$ for ease of notation. \textcolor{red}{(Lieber nicht; finde ich sehr
%verwirrend.)}
The following well-known  characterization of the Markov property for Gaussian processes
is due to Doob~\cite[Theorem~V.8.1]{Do53}. The multi-dimensional version we require
is an easy modification of the one-dimensional case.
\begin{lemma}
\label{lemma1}
Let $X$ be a centered $V$-valued Gaussian process such that $R(t,t)$ is invertible for all $t>0$.
The process $X$ satisfies the Markov property if and only if $R$ satisfies the functional equation
\begin{align}
    \label{funcequation}
    R(s,t)R(t,t)^{-1}R(t,u)=R(s,u), \quad 0\le s\le t \le u.
\end{align}
\end{lemma}
\begin{proof}
  It suffices to consider $V=\mathbb{R}^d$ with the standard inner product (see the proof of
  Proposition~\ref{prop:indep}). We follow the proof for $d=1$ given in Proposition~14.7
  of~\cite{Ka21}. The matrix $A:=R(t,t)^{-1}R(t,u)$ satisfies
  $\ee[X_t(X_u-X_tA)^\tr]=0$. The rest of the proof is the same as in~\cite{Ka21}, with~$a$
  replaced by~$A$.
\end{proof}
If~$X$ is an SSGM process, then  Lemma~\ref{le:identity}
implies that $\tilde{R}(t,t)=t^{2H}\operatorname{id}$
(recall the notation~$\tilde{R}$ from Definition~\ref{def:ssgm}).
Applying this to~\eqref{funcequation}, we obtain
\begin{align*}
t^{-2H}\tilde{R}(s,t)\tilde{R}(t,u)=\tilde{R}(s,u),
\end{align*}
and so
\[
  \tilde{R}(s/t,1)\tilde{R}(t/u,1)=\tilde{R}(s/u,1).
\]
Substituting $x:=\log(t/s)$ and $y:=\log(u/t)$ and setting $g(x):=\tilde{R}(e^{-x},1)$ yields
\begin{equation}\label{eq:sg}
g(x)g(y)=g(x+y), \quad x,y\geq 0,
\end{equation}
and hence $g$ is a one-parameter semigroup.
\begin{definition}\label{def:sg gen}
A map $g:\mathbb{R}_{\geq0} \to L(V)$ is a semigroup if $g(0)=\operatorname{id}$ and
$g(x+y)=g(x)g(y)$ for all $x,y \ge 0$.
\end{definition}
\begin{definition}\label{def:sg}
For any SSGM process $X$,  we call the one-parameter semigroup $(g_X(x))_{x \ge 0}$, where $g_X(x):=\tilde{R}(\exp(-x),1)$, its associated semigroup. If there is no ambiguity, we just write~$g$
instead of~$g_X$.
\end{definition}
We stress that our notion of associated semigroup, defined by the scaled covariance function,
is \emph{not} the transition semigroup of the Markov process, which will not be used
in this paper.
The following theorem, proven in Section~\ref{se:cov}, is our first main result.
\begin{theorem}
\label{group}
Let $H>0$,
let $(g(x))_{x \ge 0}$ be a semigroup in $L(V)$, and $\langle \cdot, \cdot \rangle$  an inner product on $V$. There exists an SSGM process~$X$
 %supported on all of~$V$
  with associated semigroup $g=g_X$ and inner product $\langle \cdot,\cdot \rangle=\langle \cdot, \cdot \rangle_X$ if and only if the contractivity estimate
\begin{equation}\label{eq:1st main thm}
  \|g(x)\|_{\mathrm{op}} \leq e^{-Hx},\quad x\geq0,
\end{equation}
holds, where $\|\cdot\|_{\mathrm{op}}$ denotes the operator norm of $\langle \cdot, \cdot \rangle$.
Furthermore, $X$ is determined uniquely in law by $g_X$ and $\langle \cdot,\cdot \rangle_X$.
\end{theorem}
\begin{definition}\label{def:wn}
We call an SSGM process white noise if its associated semigroup satisfies $g(x)=0$ for all $x > 0$.
\end{definition}
By definition of~$g$, the covariance kernel of a white noise process satisfies
$\tilde{R}(s,t)=0$ for $0<s<t$.
\begin{definition}\label{def:irr}
Let~$X$ be an SSGM process on $V$. We say that~$X$ is irreducible if for any decomposition $X=X^1+X^2$ such that~$X^1$ and~$X^2$ are independent SSGM processes supported on $V_1$ and $V_2$ respectively with $V_1 \perp V_2$, the spaces $V_1$ and $V_2$ are trivial,  i.e.\ either the entire space or $\{ 0\}$.
The orthogonality of the subspaces is with respect to $\langle \cdot,\cdot\rangle_X$.
\end{definition}
\begin{remark}
  As  $\langle \cdot,\cdot\rangle_X$ depends on~$X$ only through its law,
  so does the condition that the joint law of $(X^1,X^2)$ must satisfy
  in the preceding definition. Therefore, assuming only equality
  in distribution, $X\sim X^1+X^2$, yields an equivalent definition.
\end{remark}
We proceed to characterize all semigroups associated with centered self-similar Gaussian Markov processes.
Consider Cauchy's functional equation
\begin{equation}\label{eq:cauchy}
    f(x)+f(y)=f(x+y), \quad f:\mathbb{R}\to \mathbb{R}.
\end{equation}
%Let $D$ be the set of all solutions to this functional equation.
 %We put an equivalence relation on the set of solutions, with $\nu \sim \eta$ for $\nu,\eta \in D$ if $\nu-\eta$ is linear i.e. $\nu(x)-\eta(x)= \lambda x$ for all $x \in \mathbb{R}$ for some $\lambda \in \mathbb{R}$.
Define the rotation matrix
\[
    Q(\theta)=\begin{pmatrix}
    \cos(\theta) & \sin(\theta)\\
    -\sin(\theta) & \cos(\theta)
    \end{pmatrix}, \quad \theta \in \mathbb R,
\]
and, for even~$k$, the block-diagonal matrix
\begin{equation}\label{eq:def rot}
    Q_{k}^{\nu}(x):=\underbrace{\begin{pmatrix}Q(\nu(x)) & & 0 \\ & \ddots & \\ 0 & & Q(\nu(x))\end{pmatrix}}_{k \times k}
\end{equation}
consisting of $k/2$ rotation matrices evaluated at some solution~$\nu$ of~\eqref{eq:cauchy}.
Now we state our second main theorem, which gives representations
for the covariance kernels of all irreducible $H$-self-similar Markov processes.
Recall the notation defined in~\eqref{eq:def ip}, Definition~\ref{def:sg}, and~\eqref{eq:def rot}.
\begin{theorem}
\label{maintheorem}
Let $X$ be an irreducible SSGM process with values in $V$. Then there exists an orthonormal basis with respect to $\langle \cdot,\cdot\rangle_X$  such that with respect 
to this basis one of the following holds (we identify linear maps and matrices if the chosen basis
is clear from the context; see Remark~\ref{rem:iso} for details):
\begin{itemize}
    \item[(i)] $g_X(x)=0$ for $x>0$ and $V$ is one-dimensional (white noise case),
    \item[(ii)] $g_X(x)=\exp(Mx)$, $x\geq0$, for some matrix $M$ (the generator), or
    \item[(iii)] The dimension~$d$ is even, and $g_X(x)=Q_d^\nu(x)\exp(Mx)$, $x\geq0$, with $\nu$ a non-continuous solution to Cauchy's functional equation~\eqref{eq:cauchy} and $M$ commuting with every $Q_d^\nu(x)$. 
\end{itemize}
In the second and third case we have the contractivity condition $\|\exp(Mx)\|_{X}\le e^{-Hx}$, $x\geq0$,
where $\|\cdot\|_X$ is the operator norm on $L(V)$ induced by $\langle \cdot,\cdot \rangle_X$.
\end{theorem}
The theorem is proven at the end of Section~\ref{se:cov}.
\begin{remark}\label{rem:iso}
Let $\Gamma : V\to \mathbb{R}^d$ be an isomorphism that maps vectors to their coordinates
w.r.t.\ to the orthonormal basis from the preceding theorem. If $\tilde A$ denotes the
matrix of the inner product $\langle \cdot, \cdot \rangle_X$ in these coordinates, i.e.\
$\langle \Gamma^{-1} w,\Gamma^{-1}z\rangle_X = w^\tr \tilde{A} z$, $w,z \in \mathbb{R}^d$,
then the identification of linear maps and representation matrices we used in the theorem means that
$\|\exp(Mx)\|_{X}=\|\exp(Mx)\|_{\tilde{A}}$, where $\|\cdot\|_{\tilde A}$ is the operator norm induced by the inner product 
$(w,z)\mapsto w^\tr \tilde{A} z$ on $\mathbb{R}^d$. Similarly, in the
following remark the expression on the left 
hand side of~\eqref{eq:rem lp} is an abbreviation for
 $(\Gamma v)^\tr (M+H\,\mathrm{id})^\tr \tilde A \Gamma v$.
 \end{remark}
\begin{remark}\label{rem:lp}
  The contractivity condition $\|\exp(Mx)\|_{X}\le e^{-Hx}$, $x\geq0$, is equivalent
  to the dissipativity of the (in general non-symmetric) matrix $M+H\, \mathrm{id}$, i.e.,
  \begin{equation}\label{eq:rem lp}
    \langle(M+H\, \mathrm{id}) v, v \rangle_X \leq 0,\quad v\in V.
  \end{equation}
  This follows from the Lumer--Phillips theorem~\cite[p.~52]{Da80}, which characterizes
  contractive semigroups.
  In the special case where $\tilde{A}=\mathrm{id}$ and
   the matrix~$M$ consists of a single Jordan block, this condition
  can be simplified further. See~\cite{BaGe24b} for details, where we also point out a relation to a classical
  inequality due to Fan--Taussky--Todd~\cite{FaTaTo55}.
\end{remark}
Using Theorem~\ref{maintheorem}, we can give the associated semigroup of any SSGM process~$X$.  By the definition of irreducibility, there exists an orthogonal decomposition $V=\bigoplus_{i=1}^n V_i$ and independent irreducible SSGM processes~$X^i$ supported on~$V_i$ respectively such that $X=\sum_{i=1}^n X^i$ and $g_X=\bigoplus_{i=1}^n g_{X^i}$. By Theorem~\ref{maintheorem}, each $g_{X^i}$ has one of the forms given in Theorem~\ref{maintheorem} and satisfies the contractivity condition $\|g_{X^i}(x)\|_{X^i}\le e^{-Hx}$.
Conversely, given semigroups $g_1, \dots, g_n$ on finite-dimensional vector spaces $V_1, \dots ,V_n$ such that each $g_i$ is of the form given in Theorem~\ref{maintheorem} and satisfies the contractivity condition, then by Theorem~\ref{group}, for any inner product $\langle \cdot, \cdot \rangle$ on $\bigoplus_{i=1}^n V_i$ such that the $V_i$ are orthogonal, there exists a unique in law SSGM process~$X$ with associated semigroup $\bigoplus_{i=1}^n g_i$.

\section{Proofs}\label{se:cov}

As in Section~\ref{se:prelim}, $V$ denotes a real $d$-dimensional vector space
with inner product $\langle \cdot,\cdot \rangle$.
For one-dimensional Gaussian processes, the following result is standard. The proof
in dimension~$d$ is a straightforward extension; see~\cite[Theorem~2]{ChFaWa23}.
\begin{theorem}
\label{existence}
Let $R:\mathbb{R}_{\ge 0}\times \mathbb{R}_{\ge 0}\to L(V)$. Then there exists a  centered Gaussian process with covariance kernel $R$ iff
\begin{itemize}
    \item[(i)] $R(s,t)=R(t,s)^*$ for $s,t\ge 0$, where the adjoint is taken with respect to
     $\langle \cdot,\cdot \rangle$,
    \item[(ii)] for all $n>0$ we have $\sum_{i,j=1}^n \langle v_i,R(t_i,t_j)v_j \rangle \ge 0$ \ for $v_1,\dots,v_n \in V$ and $t_1,\dots,t_n \in \mathbb{R}_{\ge 0}$.
\end{itemize}
This process is unique in law.
\end{theorem}
The following characterization of~(ii) in Theorem~\ref{existence} seems to be new.
\begin{lemma}
\label{lemma2}
Let $R:\mathbb{R}_{\geq0}\times\mathbb{R}_{\geq0} \to L(V)$ be a kernel which satisfies \eqref{funcequation}, $R(t,t)$ being symmetric and positive definite for all $t \ge 0$ and $R(s,t)^*=R(t,s)$ for all $s,t \ge 0$. Then $R$ satisfying condition~(ii) in Theorem~\ref{existence} is equivalent to
\[
R(s,s)-R(s,t)R(t,t)^{-1}R(t,s)
\]
being positive semidefinite for all $s,t \in \mathbb{R}_{\geq0}$.
\end{lemma}
\begin{proof} We first show necessity.
Clearly it suffices to consider the problem in $\mathbb{R}^n$ equipped with the standard inner product, by applying the isomorphism $\Gamma$ from Remark~\ref{rem:iso}.
Notice that condition~(ii) in Theorem~\ref{existence} is equivalent to
\begin{equation}\label{eq:Rn}
R^n=
\begin{pmatrix}
R(t_1,t_1) & \dots & R(t_1,t_n)\\
\vdots & \ddots & \\
R(t_n,t_1) & & R(t_n,t_n)
\end{pmatrix}
\end{equation}
being positive semi-definite for all $t_1,\dots,t_n \ge 0$ and $n > 0$. For every $s,t \ge 0$ denote by $A_{s,t}$ a solution of
\[
A^\tr_{s,t}A_{s,t}=R(t,t)-R(t,s)R(s,s)^{-1}R(s,t).
\]
Let $M^n(L(V))$ the algebra of $n \times n$ matrices over the space $L(V)$. Clearly we have $R^n \in M^n$. Notice that each $M \in M^n(L(V))$ defines a linear map on $V^{n}$. Denote by $M^\tr \in M^n(L(V))$ the transpose of the operator $M$ on $V^n$.
We show by induction that the covariance operator~\eqref{eq:Rn}
%\[R^n=
%\begin{pmatrix}
%R(t_1,t_1) & \dots & R(t_1,t_n)\\
%\vdots & \ddots & \\
%R(t_n,t_1) & & R(t_n,t_n)
%\end{pmatrix}
%\]
is a positive semidefinite operator on $V^n$ for any $n \in \mathbb{N}$ and any increasing sequence $t_1,\dots ,t_n \in \mathbb{R}_{\ge 0}$. By assumption, $R^1$ is positive semidefinite. Assume that $R^n$ is positive semidefinite.
We define
\[
R^{n,n+1}:=\begin{pmatrix}
R^n & 0_n^\tr\\
0_n & \operatorname{id}
\end{pmatrix},
\]
where $0_n=(0,\dots,0)$ with $0:V \to V$ mapping everything to zero. It follows that $R^{n,n+1}$ is positive semidefinite on $\mathbb{R}^{d(n+1)}$. Let
\[
B_{n+1}:=\begin{pmatrix}
\operatorname{id}_n & v_{n+1}^\tr\\
0 & A_{t_n,t_{n+1}}
\end{pmatrix},
\]
where $\operatorname{id}_n$ is the identity in $M^n(L(V))$ and
\[
v_{n+1}:=\big(0,\dots,0,R(t_n,t_n)^{-1}R(t_n,t_{n+1})\big)\in L(V)^{1 \times n}.
\]
Then we have
\[
B_{n+1}^\tr R^{n,n+1}B_{n+1}=R^{n+1}.
\]
This shows that $R^{n+1}$ is positive semidefinite, finishing the induction step.
Conversely, by using condition~(ii) in Theorem~\ref{existence} for $n=2$, we have
\[
\langle v_1,R(s,s)v_1\rangle - 2\langle v_1,R(s,t)v_2\rangle+\langle v_2,R(t,t)v_2\rangle \ge 0.
\]
Setting $v_2 = R(t,t)^{-1}R(t,s)v_1$ yields
\[
\langle v_1, R(s,s)v_1-R(s,t)R(t,t)^{-1}R(t,s)v_1 \rangle \ge 0,
\]
which shows that $R(s,s)-R(s,t)R(t,t)^{-1}R(t,s)$ is positive semidefinite.
\end{proof}
Combining Theorem~\ref{existence}, Lemma~\ref{lemma1} and Lemma~\ref{lemma2} yields
the following result. The new aspect of Theorem~\ref{markov} is that~(iii), together with the other conditions, implies that~$R$ is a covariance kernel.
\begin{theorem}
\label{markov}
Let $R:\mathbb{R}_{\geq 0}\times \mathbb{R}_{\geq 0} \to L(V)$ be 
 such that $R(t,t)$ is
invertible for all $t>0$. Then there exists a centered Gaussian Markov process with covariance kernel $R$ if and only if
\begin{itemize}
    \item[(i)] $R(t,t)$ is positive definite for all $t>0$,
    \item[(ii)] $R(s,t)=R(t,s)^*$ for any $s,t \ge 0$, where the adjoint is taken with respect to
     $\langle \cdot,\cdot \rangle$,
    \item[(iii)] $R(s,s)-R(s,t)R(t,t)^{-1}R(t,s)$ is positive semidefinite for all $s,t \ge 0$,
    \item[(iv)] $R(s,t)R(t,t)^{-1}R(t,u)=R(s,u)$ for all $0\le s<t<u$ (Markov property).
\end{itemize}
\end{theorem}
%
%\begin{proof}[Proof of Theorem~\ref{similar}]
%Let $X$ be an $H$-self-similar Gaussian process. Since for any $a>0$ we have $(a^HX_s)\sim (X_{as})$, by uniqueness they have the same covariance kernel. This yields
%\begin{align*}
%    \langle v,R(as,at)u\rangle&= \ee[\langle v,X_{as}\rangle\langle u,X_{at}\rangle]= \ee[\langle v,a^HX_s\rangle\langle u,a^HX_t\rangle]\\
%    &=a^{2H}\ee[\langle v,X_s\rangle\langle u,X_t\rangle]= \langle v,a^{2H}R(s,t)u\rangle.
%\end{align*}
%Conversely if $R(as,at)=a^{2H}R(s,t)$, then the processes $(a^HX_s)$ and $(X_{as})$ have the same covariance %kernel and hence, again by uniqueness, are equal in law.
%\end{proof}
%
The following theorem motivates Definition~\ref{def:irr},
as it tells us that the associated semigroup of any reducible SSGM process is itself reducible in an appropriate sense.
\begin{theorem}
\label{irreducibility}
An SSGM process is irreducible if and only if there does not exist a nontrivial subspace $W\subset V$ such that $g_X(x)W= W$ and $g_X(x)W^\perp= W^\perp$ for all $x\ge 0$.
\end{theorem}
\begin{proof}
Assume first that~$X$ is \emph{not} irreducible, i.e.\ $X = X^1+X^2$ and $V=W \oplus W^\perp$, with $X^1$ supported on $W$ and $X^2$ supported on $W^\perp$. By assumption $X^1$ is independent of $X^2$. Denote by $p_1$ the orthogonal projection onto $W$ and $p_2$ the orthogonal projection onto $W^\perp$. Since $R_{X^i}(1,1)=p_iR_X(1,1)p_i=p_i$ which is the identity on $W$ and $W^\perp$ respectively, we have that $\langle \cdot,\cdot \rangle_{X^1}$ and $\langle \cdot,\cdot \rangle_{X^2}$ are just the restrictions of $\langle \cdot,\cdot \rangle_X$ to $W\times W$ and $W^\perp \times W^\perp$ respectively. Denote by $g_{X^1}(x):W \to W$ and $g_{X^2}(x):W^\perp \to W^\perp$ the associated semigroups of $X^1$ and $X^2$ respectively. We show that $g_X(x)=g_{X^1}(x)p_1+g_{X^2}(x)p_2$. Clearly, for any $v \in V$ we have $v = p_1v+p_2v$ and hence also $p_iX=X^i$. Furthermore $p_1^*=p_1$ and $p_2^*=p_2$.
 For any $v,u \in V$ we obtain
\begin{align*}
    \langle v, g_X(x)u\rangle_X &=\langle v, \tilde{R}_X(e^{-x},1)u\rangle_X=\ee[\langle v, X_{e^{-x}}\rangle_X \langle u, X_1\rangle_X]\\
    &=\ee\big[\langle v, X^1_{e^{-x}}+X^2_{e^{-x}}\rangle_X \langle u, X^1_1+X^2_1\rangle_X\big].
\end{align*}
Since $X^1$ is independent of $X^2$ we have 
\[
\ee[\langle v, X^1_{e^{-x}}\rangle_X \langle u, X^2_1\rangle_X]=\ee[\langle v, X^2_{e^{-x}}\rangle_X \langle u, X^1_1\rangle_X] = 0,
\]
and hence
\begin{align*}
\ee\big[\langle v, X^1_{e^{-x}}+X^2_{e^{-x}}\rangle_X \langle u, X^1_1+X^1_2\rangle_X\big]&=\ee\big[\langle v, X^1_{e^{-x}}\rangle_X \langle u, X^1_1\rangle_X+\langle v, X^2_{e^{-x}}\rangle_X \langle u, X^2_1\rangle_X\big]\\
&=\ee\big[\langle p_1v, X^1_{e^{-x}}\rangle_X \langle p_1u, X^1_1\rangle_X+\langle p_2v, X^2_{e^{-x}}\rangle_X \langle p_2u, X^2_1\rangle_X\big]\\
&=\langle p_1v, g_{X^1}(x)p_1u\rangle_X + \langle p_2v, g_{X^2}(x)p_2u\rangle_X\\&=\langle v, (g_{X^1}(x)p_1+g_{X^2}(x)p_2)u\rangle_X.
\end{align*}
Conversely, assume that there exists a nontrivial subspace $W$ such that $W$ and $W^\perp$ are invariant under $g_X(x)$. We show that $X^1:=p_1X$ and $X^2:=p_2X$ are independent SSGM processes. Clearly, $X =X^1+X^2$, and since $p_1,p_2$ are both linear, $X^1$ and $X^2$ are $H$-self-similar Gaussian. It remains to show that they are independent and Markovian. For Gaussian processes it is enough to show pointwise independence. For $0 \le s \le t$ and $v,u \in V$ we have
\begin{align*}
    \ee[\langle v,X^1_s \rangle_X \langle u,X^2_t\rangle_X]&=\ee[\langle p_1v,X_s \rangle_X \langle p_2u,X_t\rangle_X]\\&=\langle p_1v,\tilde{R}_X(s,t)p_2u\rangle_X=t^2H\langle p_1v,g_X(\log(t/s))p_2u\rangle_X = 0,
\end{align*}
where the last equality follows from $p_1v \in W$ and $g_X(\log(t/s))p_2u \in W^\perp$. In the same manner we obtain independence for $s>t$. For any $v,u \in V$ we have
\[
\ee[\langle v, X^1_{e^{-x}}\rangle_X \langle u, X^1_1\rangle_X]=\ee[\langle p_1v, X_{e^{-x}}\rangle_X \langle p_1u, X_1\rangle_X]=\langle v,p_1g_X(x)p_1u\rangle_X.
\]
Hence $g_{X^1}(x)=p_1g_X(x)p_1$. Analogously, we get $g_{X^2}(x)=p_2g_X(x)p_2$. It follows that $g_{X^1}(x)$ and $g_{X^2}(x)$ are semigroups, and hence $X^1$ and $X^2$ satisfy condition~(iv) in Theorem~\ref{markov}.
\end{proof}
We can now prove our main results.
\begin{proof}[Proof of Theorem~\ref{group}]
Sufficiency:
Define
\[
\tilde{R}(s,t)=t^{2H}g(\log(t/s)), \quad 0<s\leq t,
\]
and set $\tilde{R}(s,t)=\tilde{R}(t,s)^*$ for $s>t$, where the adjoint is with respect to $\langle \cdot,\cdot \rangle$. We now apply Theorem~\ref{markov} to~$\tilde R$, where~$V$ is equipped with
the inner product $\langle \cdot,\cdot \rangle$. Assumption~(i) follows from Definition~\ref{def:sg gen},
as does~(iv), using also the calculation above~\eqref{eq:sg}. Assumption~(ii) is clear
by definition of~$\tilde{R}$, and~(iii) easily reduces to
\begin{equation}\label{eq:should be pd}
\|g(x)\|_{\mathrm{op}} \le e^{-Hx}, \quad x\geq0.
\end{equation}
Hence there exists a Gaussian Markov process $X$ with covariance kernel~$\tilde{R}$.
 The resulting process~$X$ is self-similar by Theorem~\ref{similar}. The uniqueness assertion follows from Theorem~\ref{existence}.

Necessity: Let~$R$ be the covariance kernel of~$X$ w.r.t.\ $\langle \cdot,\cdot \rangle$. As
in Definition~\ref{def:sg}, we put $g_X(x)=\tilde{R}(e^{-x},1)$. 
It is easy to see that part~(iii) of Theorem~\ref{markov}
shows that~\eqref{eq:1st main thm} holds for
$\langle \cdot,\cdot \rangle_X = \langle \cdot,R(1,1)^{-1}\cdot \rangle$.
\end{proof}
The following characterization of locally bounded matrix semigroups is proven
in the companion paper~\cite{BaGe24}. It holds for~$V$ as at the beginning of Section~\ref{se:prelim},
with an arbitrary inner product $\langle\cdot,\cdot \rangle$, and we write $\|\cdot\|$
for the corresponding operator norm. Below, the theorem will be applied
to $\langle\cdot,\cdot \rangle=\langle\cdot,\cdot \rangle_X$.
\begin{theorem}
\label{semigroupdecomposition}
Let $(g(x))_{x\ge 0}$ be a semigroup in $L(V)$ such that $\|g(x)\|\le f(x)$ for all $x\ge 0$ with $f$ being locally bounded and right-continuous at 0 with $f(0)=1$. Then there exists an orthogonal decomposition $V=\bigoplus_{i=1}^n V_i$ such that each $V_i$ is invariant under $g(x)$ and either $g(x)$ is elementary  on $V_i$, or $g(x)|_{V_i} = 0$ for $x>0$. Here, \emph{elementary} means that the semigroup is of the form~(ii) or~(iii)
in Theorem~\ref{maintheorem}.
\end{theorem}
\begin{proof}[Proof of Theorem~\ref{maintheorem}]
Let $X$ be an irreducible SSGM process. By Theorem~\ref{irreducibility}, $V$ does not have a decomposition into nontrivial orthogonal subspaces which are invariant under~$g_X$. By Theorem~\ref{group}, we have $\|g_X(x)\|_{X}\le e^{-Hx}$, and hence we can apply
Theorem~\ref{semigroupdecomposition}. Because of the irreducibility this decomposition has to be trivial, and hence either $g_X(x)$ is elementary on $V$ or $g_X(x) \equiv 0$. If $g_X(x)\equiv 0$ for $x>0$ then $X$ is white noise. In
this case, for any decomposition $V=W \oplus W^\perp$ the projections $p_W(X)$ and $p_{W^\perp}(X)$ are again $H$-self-similar white noise. Since they are white noise they are also Markovian. In order for~$X$ to be irreducible, any such decomposition has to be trivial, and hence~$V$ must be one-dimensional.
\end{proof}

\section{Volterra representation}\label{se:volt}

We now assume $V=\mathbb{R}^d$, equipped with the standard inner product, and our goal is to
 give a Volterra representation of a continuous centered SSGM process~$X$ with
 semigroup $(\exp(Mx))_{x\geq0}$ satisfying
\begin{equation}\label{eq:strict lp}
  \|e^{Hx}\exp(Mx)\|_X<1, \quad x>0.
\end{equation}
Thus, we consider the non-degenerate case~(ii) of Theorem~\ref{maintheorem}, and the bound
at the end of that theorem is assumed to be strict.
% For computational reasons we want to work under the standard inner product, $\langle v,u \rangle = v^\tr u$.
We will show that there is a (deterministic) kernel $K(s,t)$ such that~$X$
equals $\int_0^\cdot K(s,\cdot)dW_s$ in distribution, where~$W$ is a Brownian motion; this
is commonly referred to as a Volterra representation of the process.
 Under the standard inner product, the covariance kernel of~$X$ takes the form
\begin{equation}\label{eq:R goal}
R(s,t)= t^{2H}\exp(M\log(t/s))A^{-1}, \quad 0<s<t,
\end{equation}
where $A:=R(1,1)^{-1}$ is a positive definite real matrix satisfying

\[
\langle u,v\rangle_X = u^\tr A v,\quad u,v \in \mathbb{R}^d.
\]
\begin{lemma}\label{le:neg def}
    Under the above assumptions, the matrix  $M^\tr A+AM+2HA$ is negative definite.
\end{lemma}
\begin{proof}
   By the Lumer--Phillips theorem (see Remark~\ref{rem:lp}), the bound
   $\|\exp((M+H\, \mathrm{id})x)\|_X\leq 1$ shows that $\langle (M+H\operatorname{id})v,v\rangle_{X}\leq 0$.
    By Theorem~4.1 in~\cite{BaGe24b}, the strict bound~\eqref{eq:strict lp} implies
    the stronger assertion
    \[
    \langle (M+H\operatorname{id})v,v\rangle_{X} < 0, \quad 0\neq v\in\mathbb{R}^d.
    \]
    (While the standard inner product is considered in~\cite{BaGe24b},
    the proof of the theorem we have just used easily extends
    to an arbitrary inner product on~$\mathbb{R}^d$.)
    We have, by symmetry of the inner product,
    \begin{align*}
    2\langle (M+H\operatorname{id})v,v\rangle_X &= \langle (M+H\operatorname{id})v,v\rangle_X+\langle v, (M+H\operatorname{id})v\rangle_X\\
    &= v^\tr(M^\tr A+AM+2HA)v. \qedhere
    \end{align*}
\end{proof}
\begin{lemma}
\label{Intequ}
Given a $d \times d$ matrix~$Q$ with $\lim_{x \to \infty} \exp(Qx)= 0$ and a $d \times d$ matrix $P$, we have, in the sense of the Cauchy principal value,
\[
P=-\int_{0}^\infty\exp(Q^\tr x)(Q^\tr P+PQ)\exp(Qx)dx.
\]
\end{lemma}
\begin{proof}
We have $\frac{d}{dx}\exp(Q^\tr x)P\exp(Qx)=\exp(Q^\tr x)(Q^\tr P+PQ)\exp(Qx)$, and hence
\[
\int_0^t  \exp(Q^\tr x)(Q^\tr P+PQ)\exp(Qx) dx = \exp(Q^\tr t)P\exp(Qt)-P.
\]
We conclude by taking the limit $t \to \infty$.
\end{proof}
The following theorem is the main result of this section.
\begin{theorem}
\label{Volterrarep}
Let $X$, $H$, $M$, and $A$ be as at the beginning of this section. In particular,
\eqref{eq:strict lp} holds. Define
\[
 K(s,t)=t^{H-1/2}A^{-1}\exp(S\log(t/s))U,\quad 0<s<t,
\]
where $U$ is a positive definite solution to $U^2=-(M^\tr A+AM+2HA)$, and $S=M^\tr +(H+1/2)\operatorname{id}$.
Then the centered Gaussian process defined by
\[
   \tilde{X}_t=\int_0^t  K(s,t)dW_s, \quad t\geq 0,
\]
where $(W_t)_{t\geq0}$ is a standard $d$-dimensional Brownian motion,
has covariance kernel~\eqref{eq:R goal}.
Therefore, it equals~$X$ in distribution.
\end{theorem}
\begin{proof}
Note that~$U$ is well-defined by Lemma~\ref{le:neg def}.
We assume $0<s\leq t$ throughout the proof. Clearly,
we have
\[
R_{\tilde{X}}(s,t)=\ee[\tilde{X}_s(\tilde{X}_t)^\tr ]=\int_0^s K(u,s)K(u,t)^\tr du,
\]
where the second equality follows from Ito's isometry. From $K(s,t)=t^{H-1/2}K(s/t,1)$ we obtain
\[
\int_0^s K(u,s)K(u,t)^\tr du=(st)^{H-1/2}\int_0^s K(u/s,1)K(u/t,1)^\tr du,
\]
and substituting $u=sr$ yields
\[
(st)^{H-1/2}\int_0^s K(u/s,1)K(u/t,1)^\tr du=s^{H+1/2}t^{H-1/2}\int_0^1 K(r,1)K(rs/t,1)^\tr dr.
\]
Now we have
\[
K(r,1)K(rs/t,1)^\tr =A^{-1}\exp(S\log(1/r))U^2\exp(S^\tr \log(1/r))\exp(S^\tr \log(t/s))A^{-1},
\]
and hence
\begin{align*}
s^{H+1/2}t^{H-1/2}&\int_0^1 K(r,1)K(rs/t,1)^\tr dr\\&=s^{H+1/2}t^{H-1/2}A^{-1}\int_{0}^1\exp(S\log(1/r))U^2\exp(S^\tr \log(1/r))dr\exp(S^\tr \log(t/s))A^{-1}\\
&=A^{-1}\left [\int_{0}^\infty\exp((M^\tr +H\operatorname{id})x)U^2\exp((M+H\operatorname{id})x)dx \right ] \ t^{2H}\exp(M\log(t/s))A^{-1},
\end{align*}
where we made the substitution $r=e^{-x}$. It suffices to calculate the integral
\[
\int_{0}^\infty\exp((M^\tr +H\operatorname{id})x)U^2\exp((M+H\operatorname{id})x)dx.
\]
Since $\|\exp((H\operatorname{id}+M)x)\|_{X}<1$ for $x>0$ by~\eqref{eq:strict lp},
we have the bound
\begin{equation}\label{eq:expo est}
  \|\exp((H\operatorname{id}+M)x)\|_{X}<(\|\exp((H\operatorname{id}+M))\|_{X})^n,
  \end{equation}
   where $n = \lfloor x \rfloor$. Clearly, this converges to 0 as $n \to \infty$.
Setting $Q := M + H\operatorname{id}$ and $P := A$ in Lemma~\ref{Intequ}, we obtain
\begin{equation}\label{eq:exp M}
\int_{0}^\infty\exp((M^\tr +H\operatorname{id})x)U^2\exp((M+H\operatorname{id})x)dx=A.
\end{equation}
Note that the integral exists, not just as a Cauchy principal value, by the
estimate~\eqref{eq:expo est}. Rewriting~\eqref{eq:exp M}
yields
\[
\int_0^s K(u,s)K(u,t)^\tr du=t^{2H}\exp(M\log(t/s))A^{-1}. \qedhere
\]
\end{proof}

\section{One-dimensional processes %in $\mathbb{R}$
}\label{se:1dim}

Specializing our multi-dimensional characterization results, Theorem~\ref{group}
and Theorem~\ref{maintheorem}, easily yields the following:
\begin{theorem}\label{thm:1dim}
If $X$ is an $\mathbb R$-valued $H$-SSGM process which is not white noise (see
 Definitions~\ref{def:ssgm} and~\ref{def:wn}), then
\[
(X_t)_{t > 0} \sim (at^{2H+M}W(t^{-2H-2M}))_{t > 0},
\]
where $M\in\mathbb R$, $a=R(1,1)^{1/2}> 0$, $H+M\le0$ and $(W(t))_{t\geq0}$ is a standard Brownian motion.
The covariance function of~$X$ is $R(s,t)=a^2t^{2H}(t/s)^M$, $0<s\leq t$.
\end{theorem}
\begin{proof}
Since $X$ is one-dimensional, it has to be irreducible. As $X$ is not white noise,
Theorem~\ref{maintheorem} yields $g_X(x)=e^{Mx}$ for some $M\in\mathbb R$. We obtain
$\tilde{R}(s,t)=t^{2H}(t/s)^M$
for $t\ge s>0$; recall that $\tilde{R}$ is the covariance kernel with respect to the inner product $\langle \cdot,\cdot \rangle_X$. By multiplying with $R(1,1)=:a^2$,  we obtain
$R(s,t)=a^2t^{2H}(t/s)^M$. It remains to check the contractivity condition on $g_X$. By 
Theorem~\ref{group} we must have $\|g_X(x)\|_{X} \le e^{-Hx}$. Since $\langle x,y \rangle_X=a^{-2}xy$, we obtain $\|g_X(x)\|_{X} = e^{Mx}$. It follows that $H+M \le 0$
(cf.\ Remark~\ref{rem:lp}). One can easily check that $at^{2H+M}W(t^{-2H-2M})$ has the same covariance kernel.
\end{proof}
The pointwise limit
 \begin{equation*}%\label{eq:wn}
    \lim_{M\to-\infty} a^2 t^{2H}(t/s)^M
    =\begin{cases}
        a^2 t^{2H} & s=t, \\
       0 & s< t,
       \end{cases}
  \end{equation*}
  of the covariance function  $R(s,t)=a^2t^{2H}(t/s)^M$
  is positive definite as well, and defines a white noise process (for $a>0$),
 which is also self-similar and  Markovian. This corresponds to case~(i) in Theorem~\ref{maintheorem},
 and we have thus characterized all one-dimensional SSGM processes. 

\begin{remark} As for Volterra representations in dimension one (cf.\ Section~\ref{se:volt}),
  note that the  borderline case $H+M=0$ leads to the degenerate process $X_t=at^H W(1)$, with $R(s,t)=a^2(st)^H$, which does not have a Volterra representation. 
  This follows from Theorem~2.2 and Remark~2.3 in~\cite{Ya15}.
  If $H+M<0$, which is~\eqref{eq:strict lp} for $d=1$, Theorem~\ref{Volterrarep}
  shows that the processes from Theorem~\ref{thm:1dim} have a Volterra representation
  with kernels
  \[
    K(s,t) = \tilde{a} t^{H-1/2}(t/s)^{M+H+1/2}, \quad 0<s<t,
  \]
  where $\tilde{a}=a\sqrt{-2(H+M)}$.
\end{remark}

A direct proof of Theorem~\ref{thm:1dim} is much easier than proving the multi-dimensional
Theorems~\ref{group} and~\ref{maintheorem},
and we thus sketch it, referring to~\cite{Ba24} for full details.
It extends the proof that fractional Brownian motion is not 
Markovian (see~\cite{Hu03} and~\cite[Theorem~2.3]{No12}).
Assume $R(1,1)=1$, and let $(g(x))_{x\geq0}$
be the semigroup from Definition~\ref{def:sg}.
By~\eqref{eq:sg}, if $x_0$ is a zero of~$g,$ then~$g$ vanishes on $[x_0,\infty).$
Suppose that there is no     interval
$[0,\varepsilon)$ with $\varepsilon>0$ on which~$g$ is positive.
Then, $g(x)=0$ for $x>0$, leading to the white noise process.
If, on the other hand, there is such an interval $[0,\varepsilon)$, 
then~\eqref{eq:sg} yields
\[
  g(x) = g(x/2^k)^{2^k}, \quad x\geq0,\ k\in\mathbb N,
\]
and by taking~$k$ sufficiently large we get $g(x)>0$ for any $x\geq 0$.
Therefore, $\varphi := \log g$ is well-defined on $[0,\infty)$,
and satisfies Cauchy's functional equation
\[
 \varphi(x+y) = \varphi(x)+ \varphi(y),  \quad x\geq y\geq 0.
\]
By symmetry of the equation, it clearly holds for $0\leq x<y$ as well.
It can be shown that the solution~$\varphi$ can be extended
to the whole real line. If~$\varphi$ is linear, then we obtain
the covariance function 
\begin{equation}\label{eq:1dim cov}
  (s,t) \mapsto (s\vee t)^{2H+M} (s\wedge t)^{-M}
\end{equation}
from Theorem~\ref{thm:1dim}. If~$\varphi$ is not linear,
then its graph is dense in $\mathbb{R}\times \mathbb{R}$ 
(see~\cite[Section~1.1]{BiGoTe87}), and it is easy
to find $0<t_1<t_2$ violating~(ii) in Theorem~\ref{existence}.
As for the condition $M\leq -H$ in Theorem~\ref{thm:1dim},
it is straightforward to check that~\eqref{eq:1dim cov}
 does not satisfy the Cauchy-Schwarz inequality for $M>-H$, hence
it cannot be a covariance function, which finishes the proof sketch.

We also mention (see~\cite{Ba24}) that there is an alternative proof of
Theorem~\ref{thm:1dim}, using the main result of~\cite{Bo82}.
Moreover, the determinant evaluation
\[
  \det \big((t_i\vee t_j)^{2H+M} (t_i\wedge t_j)^{-M}\big)_{1\leq i,j\leq k}
  =t_k^{2H} \prod_{i=1}^{k-1} t_i^{-2M}\big(t_i^{2H+2M}-t_{i+1}^{2H+2M}\big),
\]
where $0<t_1<\dots<t_k$, $k\in\mathbb N$, which follows from~\cite{Li69},
yields a non-probabilistic proof that~\eqref{eq:1dim cov}
is a positive definite kernel for $M\leq -H$.

\section{Proving non-Markovianity of one-dimensional processes}\label{se:non-markov}

As hinted at in the introduction, the initial goal of this paper was a unification
and simplification
of non-Markovianity proofs found in the literature.
 The following corollary should be applicable to
virtually any self-similar non-Markovian Gaussian
process with explicit covariance function.
 For example, it easily shows that the processes defined in \cite{BoGoTa07,RuTu06,Sg13,Sg14,Zi17}
 are not Markovian. Note that in~\cite{BoGoTa07} and~\cite{RuTu06} it is just stated that the standard
 condition~\eqref{funcequation} for Markovianity of Gaussian processes is easy to refute
 for the processes they consider.
 In~\cite{Sg13,Sg14}, the Markov property is not mentioned. In~\cite{Zi17},
 \eqref{funcequation} is used directly to prove non-Markovianity.
 \begin{definition}
   Two real functions~$f$ and~$g$ are asymptotically equivalent at
   infinity, denoted
   by $f\sim g$, if $\lim_{x\to\infty}f(x)/g(x)=1$.
 \end{definition}
\begin{corollary}\label{cor:1/2}
  Let $X$ be a  centered self-similar one-dimensional
  Gaussian process with covariance $R(s,t)=\mathbb{E}[X_s X_t]$ and $R(1,1)=1$.
  % such that  the covariance $R(t^{1/2},t)$ is not identically zero. 
  If the covariance  $R(t^{1/2},t)$ has an asymptotic expansion at infinity
  which
  contains more than one power of~$t$, or a term that 
  is not asymptotically
  equivalent to a power of~$t$ with coefficient~$1,$ then~$X$ is not Markovian.
\end{corollary}
\begin{proof}
  If $X$ is Markovian, then its covariance function satisfies
  $R(t^{1/2},t)=t^{2H+M/2}$, by Theorem~\ref{thm:1dim}. This implies the assertion, by uniqeness of asymptotic expansions (see Section~1.5
  in~\cite{deBr58}, or any other treatise on asymptotic expansions).
\end{proof}
Other functions of~$t$ than $t^{1/2}$ could be used, of course, but in all examples
from the literature we are aware of, putting $s=t^{1/2}$ seems convenient.

\begin{example}\label{ex:rl}
   The Riemann-Liouville process $X_t=\sqrt{2H}\int_0^t (t-s)^{H-1/2}dW_s$,
   where~$W$ is a Brownian motion and $H>0$, is $H$-self-similar. We have $R(s,t)=s^{2H} l((t-s)/s)$
   for $0<s\leq t$, where  (see, e.g.,~\cite{Sk19})
      \[
       l(u) := 2H \int_0^1\big((v+u)v\big)^{H-1/2}dv
       =\frac{4H u^{H-1/2}}{2H+1}
       {}_2F_1 \bigg( {{\tfrac12-H,H+\tfrac12} \atop {H+\tfrac32}} \Big| {-\frac{1}{u}}\bigg).
   \]
   Since the hypergeometric function ${}_2F_1$ is analytic at zero and ${}_2F_1(0)=1$,
   we have
   \begin{equation*}%\label{eq:l rl}
      R(t^{1/2},t) \sim \frac{4H}{2H+1} t^{3H/2-1/4},\quad t\uparrow\infty.
   \end{equation*}
   By Corollary~\ref{cor:1/2}, the process can only be Markovian if $4H/(2H+1)=1,$ i.e.\ $H=\tfrac12$.
   We conclude
   that the Riemann-Liouville process is not Markovian
   for $H\in(0,\infty)\setminus \{\tfrac12\}.$ This has been sometimes mentioned in the literature without proof.
   Another proof was recently given in Lemma~4.1 of~\cite{Wi25}.   
   Consequently, the volatility process in 
   the rough Bergomi model~\cite{BaFrGa16} from mathematical finance,
   which is defined as an exponential of the Riemann-Liouville process,
   is not Markovian.
\end{example}

The generalized fractional Brownian motion of Pang and
Taqqu~\cite{IcPaTa22,PaTa19} 
is a centered Gaussian process with two parameters
$\gamma\in[0,1)$ and $\alpha\in(-\tfrac12-\tfrac12 \gamma,
\tfrac12+\tfrac12 \gamma).$ Its covariance function is
\begin{align}
  R(s,t)&=c^2 \int_0^s(t-u)^\alpha (s-u)^\alpha u^{-\gamma}du \notag \\
  & + c^2 \int_0^\infty \big((t+u)^\alpha-u^\alpha\big)
      \big((s+u)^\alpha-u^\alpha\big)u^{-\gamma}du  \notag\\
      &=: c^2 I_1(s,t)+ c^2 I_2(s,t),\label{eq:I1 I2}
      \quad 0\leq s\leq t,
\end{align}
where $c>0$ is a normalization constant that can be expressed
in terms of the gamma and beta functions. This process is self-similar
with exponent
\[
  H = \alpha -\frac{\gamma}{2} + \frac12.
\]
The special case $\gamma=\alpha=0$ is the standard Brownian motion. The Markov
property has not been considered yet for generalized fractional Brownian motion.
\begin{theorem}
  Unless $\gamma=\alpha=0,$ the generalized fractional Brownian motion
  defined by the covariance function~\eqref{eq:I1 I2}
  is not a Markov process.
\end{theorem}
\begin{proof}
   We use Corollary~\ref{cor:1/2}, which requires the expansion
   of $R(\sqrt{t},t)$ as $t\uparrow\infty$. The integral~$I_1$ can
   be expressed by the hypergeometric function (see
   15.6.1 in~\cite{dlmf}),
   \begin{equation}\label{eq:I1 hg}
     I_1(\sqrt{t},t)
     = \frac{\Gamma(\alpha+1)\Gamma(1-\gamma)}{\Gamma(\alpha-\gamma+2)}t^{3 \alpha/2- \gamma/2 + 1/2}\,
     {}_2F_1\Big( {{-\alpha,1-\gamma} \atop {2+\alpha-\gamma}} \Big| \frac{1}{\sqrt{t}}\Big).
   \end{equation}
   In the integral~$I_2$, we substitute $u=tv$ and
   split the integration range,
   \begin{align}
     I_2(\sqrt{t},t)&=t^{2\alpha-\gamma+1}\Big(\int_0^{t^{-1/4}}
     +\int_{t^{-1/4}}^\infty\Big)
     \big((1+v)^\alpha-v^\alpha\big) \big((t^{-1/2}+v)^\alpha-v^\alpha\big)
       v^{-\gamma}dv \notag \\
       &=: t^{2\alpha-\gamma+1}\big(I_{21}(\sqrt{t},t)+I_{22}(\sqrt{t},t)\big).
       \label{eq:I split}
   \end{align}
   For $\alpha\geq0$ and $0\leq v\leq t^{-1/4},$ we have
   $(1+v)^\alpha-v^\alpha\sim 1,$ and thus 
   \[
     I_{21}(\sqrt{t},t)\sim
     \int_0^{t^{-1/4}}
       \big((t^{-1/2}+v)^\alpha-v^\alpha\big)
       v^{-\gamma}dv,
   \]
   which yields
   \begin{equation}\label{eq:I21}
      %t^{2\alpha-\gamma+1}
      I_{21}(\sqrt{t},t)\sim
      \frac{1}{1-\gamma}t^{-\alpha/2+\gamma/4-1/4}\,
        {}_2F_1\Big( {{-\alpha,1-\gamma} \atop {2-\gamma}} \Big| -t^{1/4}\Big)
        - \frac{t^{-\alpha/4+\gamma/4-1/4}}{\alpha-\gamma+1}.
   \end{equation}
   If $\alpha<0,$ then $(1+v)^\alpha-v^\alpha\sim -v^\alpha,$
   and we obtain
   \begin{align*}
       I_{21}(\sqrt{t},t)&\sim
     \int_0^{t^{-1/4}}
       \big(v^\alpha-(t^{-1/2}+v)^\alpha\big)
       v^{\alpha-\gamma}dv\\
       &=\frac{t^{-\alpha/2+\gamma/4-1/4}}{2\alpha-\gamma+1}
       -\frac{t^{-3\alpha/4+\gamma/4-1/4}}{\alpha-\gamma+1}
        {}_2F_1\Big( {{-\alpha,\alpha-\gamma+1} \atop {\alpha-\gamma+2}} \Big| -t^{1/4}\Big).
   \end{align*}   
   For $t^{-1/4}<v,$ we have $v\sqrt{t}\uparrow\infty,$ hence
   \begin{equation}\label{eq:binom}
     (t^{-1/2}+v)^\alpha-v^\alpha = \frac{\alpha}{\sqrt{t}}v^{\alpha-1}
     +\frac{\alpha(\alpha-1)v^{\alpha-2}}{2t}
     +O(v^{\alpha-3}t^{-3/2}),
   \end{equation}
   and thus
   \begin{equation}\label{eq:I22}
      I_{22}(\sqrt{t},t)\sim
      \frac{\alpha}{\sqrt{t}} \int_{t^{-1/4}}^\infty\big((1+v)^\alpha
      -v^{\alpha}\big) v^{\alpha-\gamma-1}dv.
   \end{equation}
   If $\alpha>\gamma$, this is $\alpha c_0 t^{-1/2}+O(t^{-3/4}),$
   because the integral
   \[
       c_0:=\int_0^\infty  \big((1+v)^\alpha-v^\alpha\big) v^{\alpha-\gamma-1}dv
   \]
   exists.
   If $\alpha\leq\gamma,$ then~\eqref{eq:I22} implies
   \begin{align}
       I_{22}(\sqrt{t},t)&\sim
       \frac{\alpha}{\sqrt{t}}\int_{t^{-1/4}}^{1/2}\big(1+O(v)
      -v^{\alpha}\big) v^{\alpha-\gamma-1}dv\\
      &\sim\frac{\alpha}{\sqrt{t}}\int_{t^{-1/4}}^{1/2}(1
      -v^{\alpha}) v^{\alpha-\gamma-1}dv \label{eq:I22 as} \\
      &=\begin{cases}
        \frac{\alpha}{\sqrt{t}}\Big(
        \frac{t^{-\alpha/2+\gamma/4}}{2\alpha-\gamma}
        +\frac{t^{-\alpha/4+\gamma/4}}{\gamma-\alpha}
        +O(1)\Big), & \alpha < \gamma, \notag \\
        \frac{\gamma \log t + O(1)}{4\sqrt{t}}, & \alpha=\gamma.
      \end{cases} \notag
     % &\sim \frac{\alpha}{\gamma-\alpha}t^{-\alpha/4+\gamma/4-1/2}. 
   \end{align}

   {}From these estimates  and
    the well-known expansion of the hypergeometric function,
   it is now a straightforward matter to check that the expansion of $R(\sqrt{t},t)$
   as $t\uparrow \infty$ has at least two terms. A case distinction according
   to the values of $\alpha$ and $\gamma$ is required, but as all cases are easy
   and rather similar, we just make the case where $\alpha>0$ and $\alpha>\gamma$
   explicit, and make this assumption throughout the rest of the proof.
      By~\eqref{eq:I1 hg} and {\S}15.12(i) of~\cite{dlmf},
   \begin{multline}\label{eq:I1 as}
      I_1(\sqrt{t},t)=
      \frac{\Gamma(\alpha+1)\Gamma(1-\gamma)}{\Gamma(\alpha-\gamma+2)}
      t^{3\alpha/2-\gamma/2+1/2}\\
      +\frac{\alpha(\gamma-1)\Gamma(\alpha+1)\Gamma(1-\gamma)}{\Gamma(\alpha-\gamma+3)}
      t^{3\alpha/2-\gamma/2}
      \big(1+o(1)\big).
   \end{multline}
   Using {\S}15.12(i) of~\cite{dlmf} in~\eqref{eq:I21} shows that the first term
   of the expansion of ${}_2F_1$ cancels with the power of~$t$
   at the end of~\eqref{eq:I21}, and the second term of the
   expansion of ${}_2F_1$ yields
   \[
     I_{21}(\sqrt{t},t) \sim \frac{\alpha}{\alpha-\gamma}
     t^{-\alpha/4+\gamma/4-1/2}.
   \]
   As noted after~\eqref{eq:I22}, $I_{22}(\sqrt{t},t)=\alpha c_0 t^{-1/2}
   +O(t^{-3/4}).$
   (Reasoning analogously to~\eqref{eq:I22 as}, it is easy
   to see that the term resulting from the second term
   on the right hand side of~\eqref{eq:binom} is of lower order.)
   Since
   \[
     -\frac{\alpha}{4} + \frac{\gamma}{4}-\frac12 > -\frac34,
   \]
   we have
   \begin{equation}\label{eq:I21 expans}
     I_{21}(\sqrt{t},t)+I_{22}(\sqrt{t},t)
     =\alpha c_0 t^{-1/2}
     +\frac{\alpha}{\alpha-\gamma}
     t^{-\alpha/4+\gamma/4-1/2}
     \big(1+o(1)\big).
   \end{equation}
   Then, by
   \[
     2\alpha-\gamma+1
      +\Big( {-\frac{\alpha}{4}} + \frac{\gamma}{4}-\frac12 \Big)
      >\frac{3\alpha}{2}-\frac{\gamma}{2}+\frac12,
   \]
   we see from~\eqref{eq:I split}, \eqref{eq:I1 as} and~\eqref{eq:I21 expans}
   that~$I_1$ does not contribute to the first two terms
   of the expansion of~\eqref{eq:I1 I2}, and so   
   \begin{equation*}
      R(\sqrt{t},t)=
     \alpha c_0 c^2 t^{2\alpha-\gamma+1/2}
     +\frac{\alpha c^2}{\alpha-\gamma}
     t^{7\alpha/4-3\gamma/4+1/2}
      \big(1+o(1)\big), \quad t\uparrow \infty.
   \end{equation*}
   As both coefficients in this expansion are clearly non-zero, and $\alpha>\gamma$
   implies that the two exponents are different, we can
   apply Corollary~\ref{cor:1/2}.
\end{proof}

\bigskip

   {\bf Acknowledgement.} 
   We are indebted to  Christa Cuchiero, Christian Krattenthaler,
   Mikhail Lifshits, and the anonymous reviewers for very helpful comments.

\end{document}